\documentclass{naw-latex}
\input epsf
\def\relabelbox{%
  \hbox\bgroup%
}%
\def\endrelabelbox{%
\egroup%
}%
\def\relabel #1#2 {%
  \special{ps:/a {} def}%
  \smash{\rlap{#2}}%
}%
\def\adjustrelabel <#1,#2> #3#4 {%
  \special{ps:/a {} def}%
  \smash{\rlap{\kern #1 \raise #2\hbox{#4}}}%
}%
\def\extralabel <#1,#2> #3 {\smash{\rlap{\kern #1 \raise #2\hbox{#3}}}}%

\usepackage{graphicx}
\usepackage{color}

\newcommand{\R}{{\mathbb R}}

\newenvironment{Fermat}%
               {\endgraf\smallbreak
                \noindent\textbf{Fermat's Principle}\sffamily}%
               {\endgraf\smallskip}

\newenvironment{Huygens}%
               {\endgraf\smallbreak
                \noindent\textbf{Huygens' Principle}\sffamily}%
               {\endgraf\smallskip}

\bibitem{andr94}{C. D. Andriesse,
{\it Titan kan niet slapen: Een biografie van Christiaan Huygens},
Uitgeverij Contact, Amsterdam, 1994.}
\bibitem{quest}{W. J. H. Andrewes, ed.,
{\it The Quest for Longitude},
Proceedings of the Longitude Symposium (Harvard University, 1993),
Collection of Historical Scientific Instruments,
Harvard University, Cambridge, MA, 1996.}
\bibitem{bern}{{\it Die Streitschriften von Jacob und Johann Bernoulli:
Variationsrechnung}, Birk\"auser Verlag, Basel, 1991.}
\bibitem{fauv00}{J. Fauvel et al., eds.,
{\it Oxford Figures: 800 Years of the
Mathematical Sciences},
Oxford University Press, Oxford, 2000.}
\bibitem{geig99/01}{H. Geiges,
Facets of the cultural history of mathematics, oratie, Universiteit
Leiden, 1999; published in
{\it European Review} {\bf 8} (2000), 487--497.}
\bibitem{geig01}{H. Geiges,
A brief history of contact geometry and topology,
{\it Expo. Math.} {\bf 19} (2001), 25--53.}
\bibitem{hawa00}{E. Hairer et G. Wanner,
{\it L'analyse au fil de l'histoire},
Springer-Verlag, Berlin, 2000.}
\bibitem{hitr96}{S. Hildebrandt and A. J. Tromba,
{\it The Parsimonious Universe: Shape and Form in the Natural World},
Springer-Verlag, New York, 1996.}
\bibitem{huyg73}{C. Huygens,
{\it Horologium oscillatorium sive de motu pendulorum ad
horologia aptato demonstrationes geometricae},
Paris, 1673; English translation:  {\it The pendulum clock or geometrical
demonstrations concerning the motion of pendula as applied to clocks},
Iowa State University Press, 1986.}
\bibitem{huyg90}{C. Huygens,
{\it Trait\'e de la lumi\`ere}, Paris, 1690; original and Dutch
translation: {\it Verhandeling over het licht}, Epsilon Uitgaven,
Utrecht, 1990.}
\bibitem{ibm}{IBM Poster {\it Bedeutende Mathematiker}.}
\bibitem{maur93}{Mauritshuis,
{\it Illustrated General Catalogue},
Meulenhoff, Amsterdam, 1993.}
\bibitem{melv51}{H. Melville,
{\it Moby Dick}, 1851.}
\bibitem{sobe98}{D. Sobel and W. J. H. Andrewes,
{\it The Illustrated Longitude},
Fourth Estate Limited, London, 1998.}
\bibitem{thur94}{W. P. Thurston,
On proof and progress in mathematics,
{\it Bull. Amer. Math. Soc. (N.S.)} {\bf 30} (1994), 161--177.}
\bibitem{yode88}{J. G. Yoder,
{\it Unrolling Time: Christiaan Huygens and the Mathematization of Nature},
Cambridge University Press, Cambridge, 1988.}

\begin{document}

\StartArtikel[Titel={Christiaan Huygens and Contact Geometry},
          AuteurA={Hansj\"org Geiges},
          AdresA={Mathematisches Institut\crlf
                  Universit\"at zu K\"oln\crlf
                  Weyertal 86--90\crlf
                  D-50931 K\"oln\crlf},
          EmailA={geiges@math.uni-koeln.de},
          kolommen={2}
          ]

\StartLeadIn
This article is based on the author's inaugural lecture at the
Universit\"at zu K\"oln on 24 January 2003.
\StopLeadIn

\vspace{2mm}

\begin{flushright}
\begin{minipage}{7cm}
{\small ``--- Oui, voil\`a le g\'eom\`etre! Et ne crois pas que les\\
g\'eom\`etres n'aient pas \`a s'occuper des femmes!''}
\end{minipage}

\vspace{2mm}

\begin{minipage}{5cm}
{\small Jean Giraudoux,\\
{\it La guerre de Troie n'aura pas lieu}}
\end{minipage}
\end{flushright}

\vspace{2mm}

\onderwerp{Introduction}
For me, the most evocative
painting in the Mauritshuis in Den Haag has always been
{\em Het meisje met de oorbel}, even before a novel and a film
turned the girl into something of a pop icon. However, that museum is
the home to another portrait that cannot fail to attract the attention
of any scientifically interested visitor, and one where the
identity of (some of) the portrayed, like
in Vermeer's famous painting, is shrouded in mystery. I am speaking of
Adriaan Hanneman's {\em Portret van Constantijn Huygens en zijn kinderen}
(Figure~1).
This family portrait depicts C.~Huygens (1596--1687) ---
``the most versatile and
the last of the true Dutch Renaissance virtuosos'' (Encyclopaedia Britannica),
whose most notable contributions lay in the fields of diplomacy
and poetry ---, together with his five children. Among them is Christiaan
Huygens (1629--1695),
who would go on to become one of the most famous mathematicical
scientists of his time, later to be characterised as
``ein Junggeselle von hervorragendem Charakter
und au{\ss}ergew\"ohnlicher Intelligenz''~\cite{ibm}.
While I expound some of the mathematical themes of Christiaan Huygens' life
and hint at their relation to modern contact geometry,
I leave the reader to ponder the question just which of the four boys
in the family portrait shows that intellectual promise, a question to which
I shall return at the end of this article.

\vspace{2mm}
\begin{minipage}{8cm}
\begin{center}
\epsfxsize 8cm \epsfbox{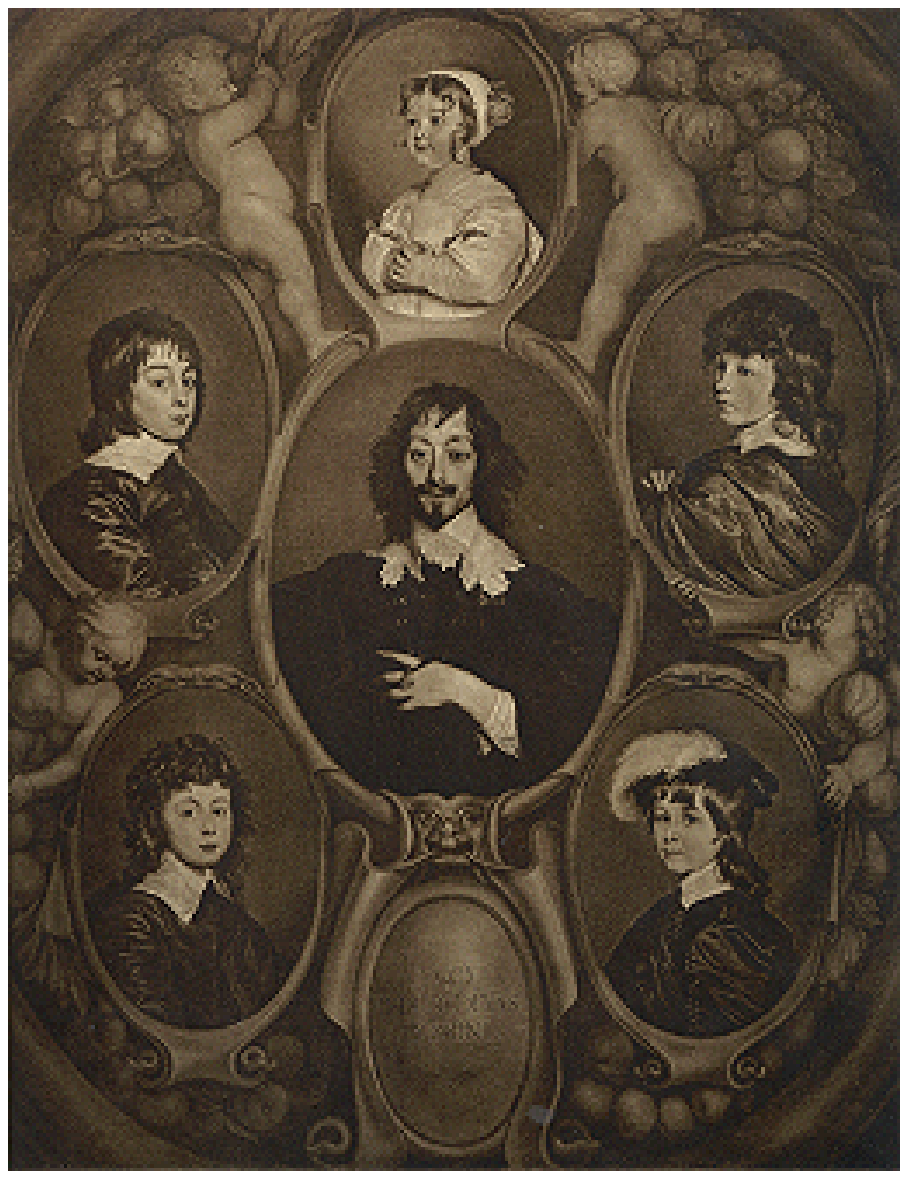}

Figure 1
\end{center}
\end{minipage}
\vspace{2mm}

An inaugural lecture is not only an opportunity to present one's
field of research to a wider public, it also allows one to reflect
on the standing of mathematics within the general intellectual discourse.
On an earlier occasion of this kind~\cite{geig99/01}
I have not been overly optimistic in this respect, and I have no reason to
qualify anything I said there. Still, it is worth
remembering that there have been even more precarious times for
mathematics. In~\cite{fauv00} we read that ``The new Savilian
professor [Baden Powell, Savilian professor of geometry at the University
of Oxford 1827--1860] was shocked and dismayed by the low esteem
accorded to mathematics in the University. He had been advised not
to give an inaugural lecture on arrival, as he would almost
surely not attract an audience.''

\vspace{1mm}

{\bf Disclaimer.} A foreigner, even one who has lived in the Netherlands
for several years,
is obviously carrying tulips to Amsterdam (or whatever the appropriate turn
of phrase might be) when writing about Christiaan Huygens in a Dutch journal.
Then again, from a visit to the Huygensmuseum Hofwijck in Voorburg near Den
Haag I gathered that in the Netherlands the fame of Constantijn Huygens
tends to outshine that of his second-eldest son. Be that as it may, this
article is intended merely as a relatively faithful record 
of my inaugural lecture (with some mathematical details added)
and entirely devoid of scholarly aspirations.
I apologise for the fact that the transcription from the spoken
to the written word has taken rather longer than anticipated.
\onderwerp{The best slide for twins}
Imagine that you are trying to connect two points $A,B$ in a vertical
plane by a slide along which a point mass $M$ will move, solely under
the influence of gravitation, in shortest time from $A$ to~$B$
(see Figure~2). This
is the famous brachistochrone problem (from Greek $\beta\rho\acute{\alpha}
\chi\iota\sigma\tau o\varsigma$ = shortest, $\chi\rho\acute{o}\nu o\varsigma$
= time), posed by Johann Bernoulli in 1696 in rather more
erudite language: ``Datis in plano verticali duobus punctis $A$ \& $B$
assignare Mobili $M$ viam AMB, per quam gravitate sua descendens \& moveri
incipiens a puncto $A$, brevissimo tempore perveniat ad alterum
punctum~$B$.'' ({\it Problema novum ad cujus solutionem Mathematici
invitantur}, Joh.\ Op.~XXX (pars), \cite[p.~212]{bern}).

\vspace{2mm}
\begin{minipage}{8cm}
\begin{center}
\centerline{\relabelbox\small
\epsfxsize 4cm \epsfbox{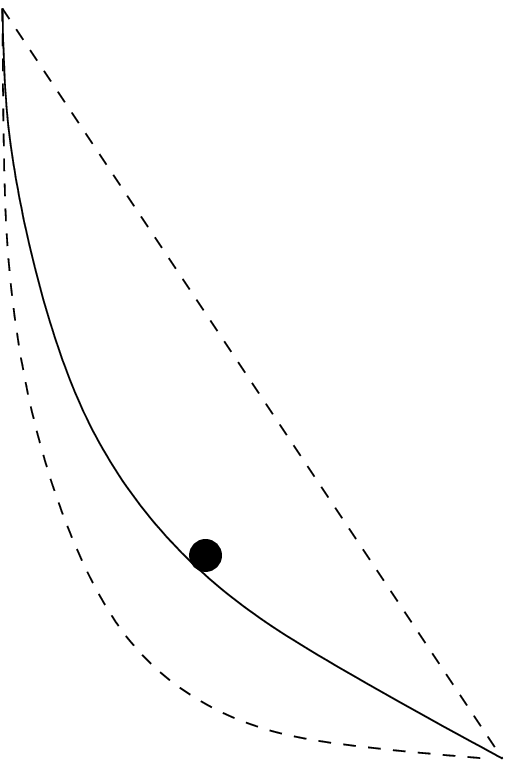}
\extralabel <-2.3cm, 1.8cm> {$M$}
\extralabel <-4.5cm, 5.8cm> {$A$}
\extralabel <0.1cm, -0.1cm> {$B$}
\endrelabelbox}

Figure 2
\end{center}
\end{minipage}
\vspace{2mm}

A related problem is to find the slide connecting the points $A$ and
$B$ in such a way that one will reach the endpoint B of the slide in the
same amount of time, no matter where on the slide one starts.
This is known as the tautochrone problem. 

Rather surprisingly, it turns out that the solution to either question
is one and the same curve, the so-called {\it cycloid}.
This is obviously the best slide a doting uncle can build for his
twin nephews: not only will their slide be faster than anybody else's;
if both of them start at the same time at any two points of the slide,
they will reach the bottom of the slide simultaneously. This gives them
the chance and the time to fight over other things.

In 1697 Jacob Bernoulli responded to the challenge set by his brother
concerning the brachistochrone with
a paper bearing the beautiful title {\it Solutio Problematum Fraternorum,
una cum Propositione reciproca aliorum},
Jac.\ Op.~LXXV~\cite[pp.~271--282]{bern}. Johann's own solution appeared
the same year (Joh.\ Op.~XXXVII, \cite[pp.~263--270]{bern}).
The tautochrone problem had been solved by Christiaan
Huygens as early as 1657, but the solution was not published until 1673
in his famous {\it Horologium Oscillatorium}~\cite{huyg73}, cf.~\cite{yode88}.

\onderwerp{The cycloid}

The cycloid is the locus traced out by a point on the rim of a circle
as that circle rolls along a straight line (Figure~3). Choose cartesian
coordinates in the plane such that the circle rolls along the $x$--axis,
with the point on the rim initially lying at the origin $(0,0)$.
Let $a$ be the radius of the circle. When the circle has turned through an
angle~$t$, its centre lies at the point $(at,a)$, and so a
parametric
description of the cycloid is given by
\[ x(t)=a(t-\sin t),\;\; y(t)=a(1-\cos t).\]

\vspace{2mm}
\begin{minipage}{8cm}
\begin{center}
\epsfxsize 8cm \epsfbox{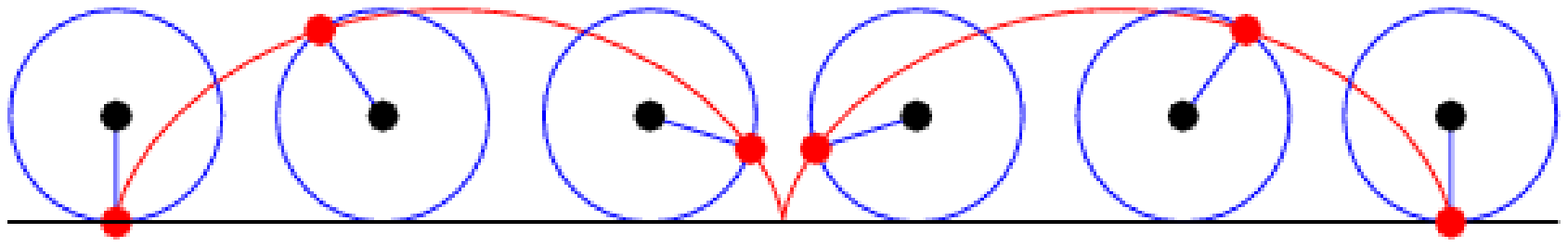}

Figure 3  (from {\tt http://mathworld.wolfram.com})
\end{center}
\end{minipage}
\vspace{2mm}

The cycloidal slide is obtained by turning this curve upside down. It is
convenient to effect this by reversing the direction of the $y$--coordinate,
while keeping the parametric equations unchanged. Given two points
$A=(0,0)$ and $B=(b_1,b_2)$ with $b_1>0$, $b_2\geq 0$ in the $xy$--plane,
there is a unique radius $a$ and angle $t_B\in (0,2\pi ]$ such that
$A=(x(0),y(0))$ and $B=(x(t_B),y(t_B))$. There are various ways to see this,
the following is due to Isaac Newton, cf.~\cite[p.~43]{bern}: Draw any
cycloid starting at~$A$, and let $Q$ be its intersection with the
straight line segment $AB$. Then expand the cycloid by a factor $AB/AQ$.
(Here and below I use the same symbol $AB$ to denote a curve or line
segment between two points $A,B$, as well as the length of that
segment, provided the meaning is clear from the context.)

For some of the reasonings below I shall assume implicitly that $t_B\leq \pi$,
so that the cycloidal segment connecting $A$ and $B$ does not have any
upward slope; this is equivalent to requiring $b_2\geq 2b_1/\pi$.

The brachistochrone and tautochrone problems were two of the most challenging
geometric questions of 17th century mathematics, attracting the
attention of the most famous (and cantankerous) mathematicians
of that time, including the Marquis de L'Hospital, Leibniz,  and
Newton. As a result, these problems were
the source of acrimonious battles over priority --- the publications of
the Bernoulli brothers on this topic have even been published in a collection
bearing the title {\it Streitschriften}~\cite{bern}. This was not the only
occasion when the cycloid was the object of desire in a mathematical
quarrel, and so this curve has often been dubbed the `Helen of Geometers'.

The following allusion to the tautochronous property of the cycloid
in Herman Melville's {\it Moby Dick}~\cite[Chapter~96, The Try-Works]{melv51}
shows that there
were happy times when the beauty of mathematics had to some degree
entered popular consciousness: ``[The try-pot\footnote{A pot for trying oil
from blubber.}] is a place also for profound mathematical
meditation. It was in the left hand try-pot of the
{\it Pequod}\footnote{Captain Ahab's ship, named after an Indian people.},
with the soapstone diligently circling around me, that I was first
indirectly struck by the remarkable fact, that in geometry all bodies
gliding along the cycloid, my soapstone for example, will descend from
any point in precisely the same time.''

\onderwerp{The cycloidal pendulum}

Besides the discovery of the true shape of Saturn's rings and one of its
moons, namely Titan, Christiaan Huygens' most important scientific
contributions are his theory of light, based on what has become known
as Huygens' principle (discussed in the next section),
and his development of a pendulum clock starting from his proof of
the tautochronous property of the cycloid.

At the time of Huygens, pendulum clocks were built (as they usually are
today) with a simple circular pendulum. The problem with such a pendulum is
that its frequency depends on the amplitude of the oscillation. With regard to
the pendulum clock in your living room this is no cause for concern,
since there the amplitude stays practically constant. But arguably the
most outstanding problem of applied mathematics at that time was to build
a clock that was also reliable in more adverse conditions, say on a ship
sailing through gale force winds. Why are such accurate clocks important?

As is wryly remarked in the introduction to the
lavishly illustrated proccedings of the Longitude Symposium~\cite{quest},
``Traveling overseas, we now complain when delayed for an hour:
we have forgotten that once there were problems finding continents''.
Indeed, how was it possible to determine your
exact position at sea (or anywhere else, for that matter), prior to the days
of satellite-based Global Positioning Systems? Mathematically the answer
is simple (at least on a sunny day): Observe when the sun reaches
its highest elevation. This will be noon local time. Moreover, the
angle $\alpha$ of elevation will give you the latitude: If the axis of
the earth's rotation were orthogonal to the plane in which the earth
moves around the sun, that latitude would simply be $90^{\circ}-\alpha$.
In order to take the tilting of the earth's axis by $23^{\circ}$
into account, one needs to
adjust this by an angle that depends on the date, varying between
$0^{\circ}$ at the equinoxes and $\pm23^{\circ}$ at the solstices.

The longitude, on the other hand, cannot be determined from this
observational data alone. Indeed, the actual value of the longitude at
any given point is a matter of convention. The fact that the zero meridian
passes through Greenwich is a consequence of the scientific achievements
and geopolitical power of the British, not astronomy.
However, if you keep a clock with
you that shows accurate Greenwich time, and you bear in mind that
the earth rotates by a full $360^{\circ}$ in 24 hours, then multiplying
the difference between your local time and that shown on the clock
by $15^{\circ}/h$ will determine your longitude relative to that of Greenwich.

All the practical problems involved in building such an accurate clock
were first solved by John Harrison in 1759, cf.~\cite{quest} and the
thrilling account of Harrison's life in~\cite{sobe98}.

From a mathematical point of view, the question addressed by Huygens concerned
the most interesting aspect of these practical problems:
Is it possible to devise a pendulum whose
frequency does not depend on the amplitude of the pendular motion?
The hardest part of this question is to find the tautochronous curve,
along which the pendulum mass should be forced to move. This Huygens
established to be the cycloid. He further observed that one could
make the pendulum move along a cycloid by restricting the swinging motion
of the pendulum between appropriately shaped plates.

\vspace{3mm}
\begin{minipage}{8cm}
\begin{center}
\centerline{\relabelbox\small
\epsfxsize 8cm \epsfbox{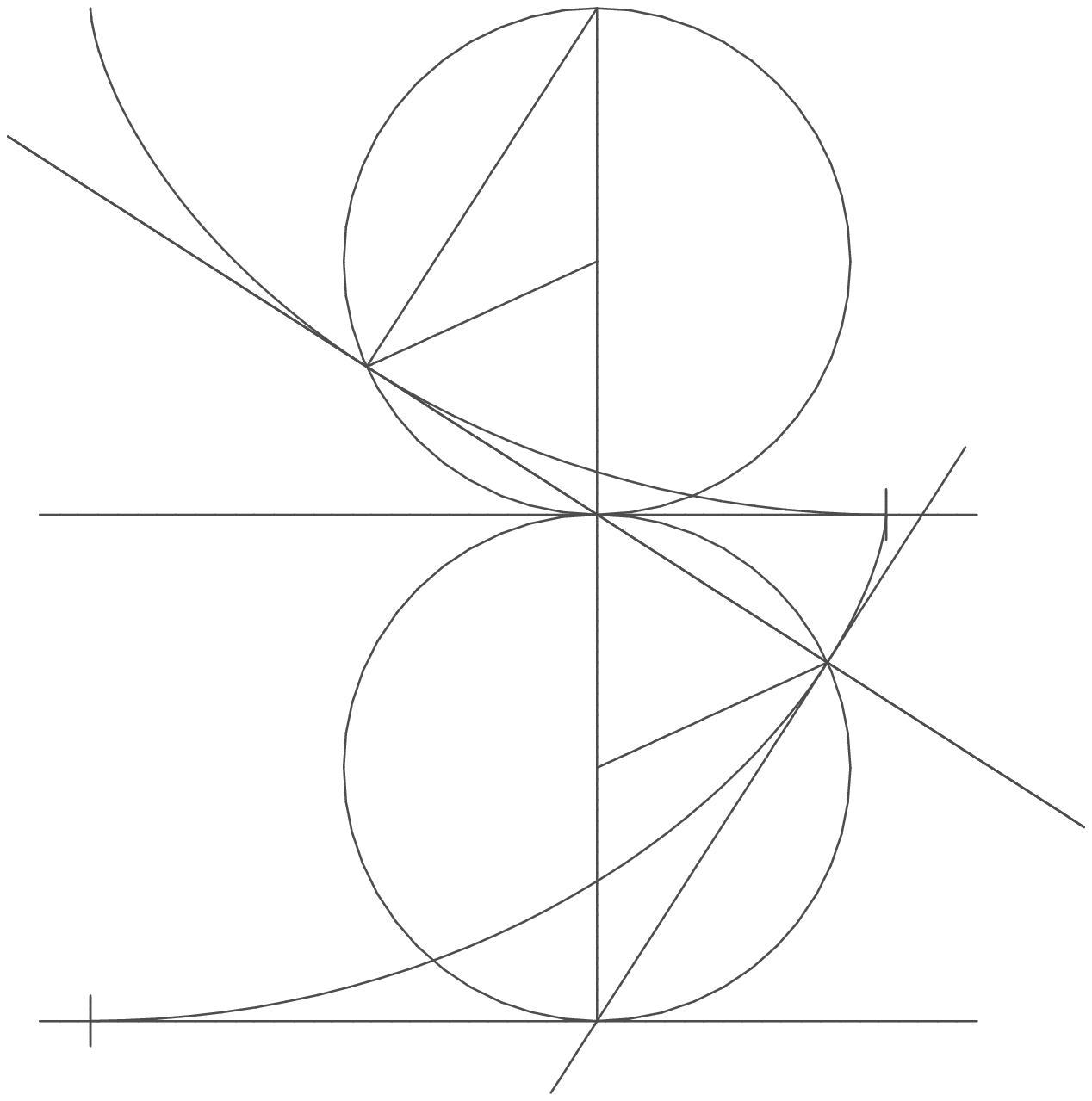}
\extralabel <-1.8cm, 3.2cm> {$M$}
\extralabel <-5.9cm, 5.1cm> {$M'$}
\extralabel <-1.7cm, 4.6cm> {$A$}
\extralabel <-7.6cm, 0.0cm> {$B$}
\extralabel <-3.7cm, 0.2cm> {$P$}
\extralabel <-4.1cm, 4.0cm> {$P'$}
\extralabel <-4.1cm, 2.4cm> {$O$}
\extralabel <-7.9cm, 8.0cm> {$B'$}
\endrelabelbox}

Figure 4
\end{center}
\end{minipage}
\vspace{2mm}

Take a look at Figure~4 (kindly provided by Manfred Lehn).
Here $AB$ is (half) a cycloid, along which the
pendulum mass~$M$, attached to the string $B'M$, is supposed to move.
This means that we require this string to be tangent to the curve $B'A$
at the point~$M'$, and the length $B'M$ to equal $B'A$,
the length of the pendulum. In other words, the
cycloid $AB$ is given by tightly unrolling
(whence the title of~\cite{yode88}) a string from the curve $B'A$.
If the pendulum is forced to swing between two plates shaped like $B'A$,
then the pendulum mass will move along the cycloid, as desired.

Such a curve $AB$ obtained by unrolling a string from a curve
$B'A$ is called the {\it involute} of~$B'A$ (and $B'A$ the
{\it evolute} of~$AB$). So the second question
faced by Huygens was: Which curve has the cycloid
as its involute? Rather miraculously, the answer is again: the cycloid.
Here is the geometric proof: Let $AB$ be the cycloid traced out by the
point $M$ as the lower circle in Figure~4 rolls to the left
along the horizontal line between the two circles (with $M=A$ at $t=0$),
and $B'A$ the cycloid traced out by the point $M'$ as the upper circle rolls
to the right along a horizontal line through $B'$ (with $M'=B'$ at $t=0$).
With the defining equations for the cycloids as in the previous section,
the situation shown in the figure corresponds to $t=t_0$ for some
$t_0\in [0,\pi ]$ for the lower circle and $t=\pi-t_0$ for the
upper circle.

The velocity (with respect to the
parameter~$t$) of the point $M$ can be split into two vector components
of length~$a$: one in horizontal direction, corresponding to the
speed of the centre of the circle, and one in the direction tangent
to the circle, corresponding to the angular speed of the rolling circle.
An elementary consideration shows that the line $MP$ bisects the
angle between these two directions, and so this line constitutes the
tangent line to the cycloid at~$M$. Analogously, the line $M'P'$ is
the tangent line to the cycloid $B'A$ at~$M'$. By symmetry
of the construction, the line $M'P'$ passes through $M$.
In order to conclude that $AB$ is
the involute of $B'A$ it therefore suffices to show that the length of
the cycloidal segment $M'A$ equals the length of the line segment $M'M$.
Also observe that, by the theorem of Thales,
the line $M'M$ is orthogonal to the tangent line
$MP$ at~$M$ ; this is a general phenomenon for an involute.

The angle $\angle MOP'$ equals~$t_0$, so the law of cosines applied to
the triangle $OMP'$ yields
\[ (P'M)^2=2a^2-2a^2\cos t_0=4a^2\sin^2\frac{t_0}{2},\]
hence
\[ M'M=2P'M=4a\sin \frac{t_0}{2}.\]

On the other hand, from the defining equations of the cycloid we have
\[ \dot{x}^2+\dot{y}^2=a^2(1-\cos t)^2+ a^2\sin^2t=4a^2\sin^2\frac{t}{2},\]
whence
\[ M'A=\int_{\pi -t_0}^\pi 2a\sin \frac{t}{2}\; dt
=4a\cos \frac{\pi -t_0}{2}= 4a\sin \frac{t_0}{2},\]
that is, $M'M=M'A$, which was to be shown.

Huygens did not stop at these theoretical considerations,
but proceeded to construct an actual pendulum
clock with cycloidal plates. The construction plan from Huygens'
{\it Horologium Oscillatorium}, with the cycloidal plates
indicated by `FIG.~II', is shown in Figure~5.
A replica of this clock can be seen in the Huygensmuseum Hofwijck.

\vspace{2mm}
\begin{minipage}{8cm}
\begin{center}
\epsfxsize 8cm \epsfbox{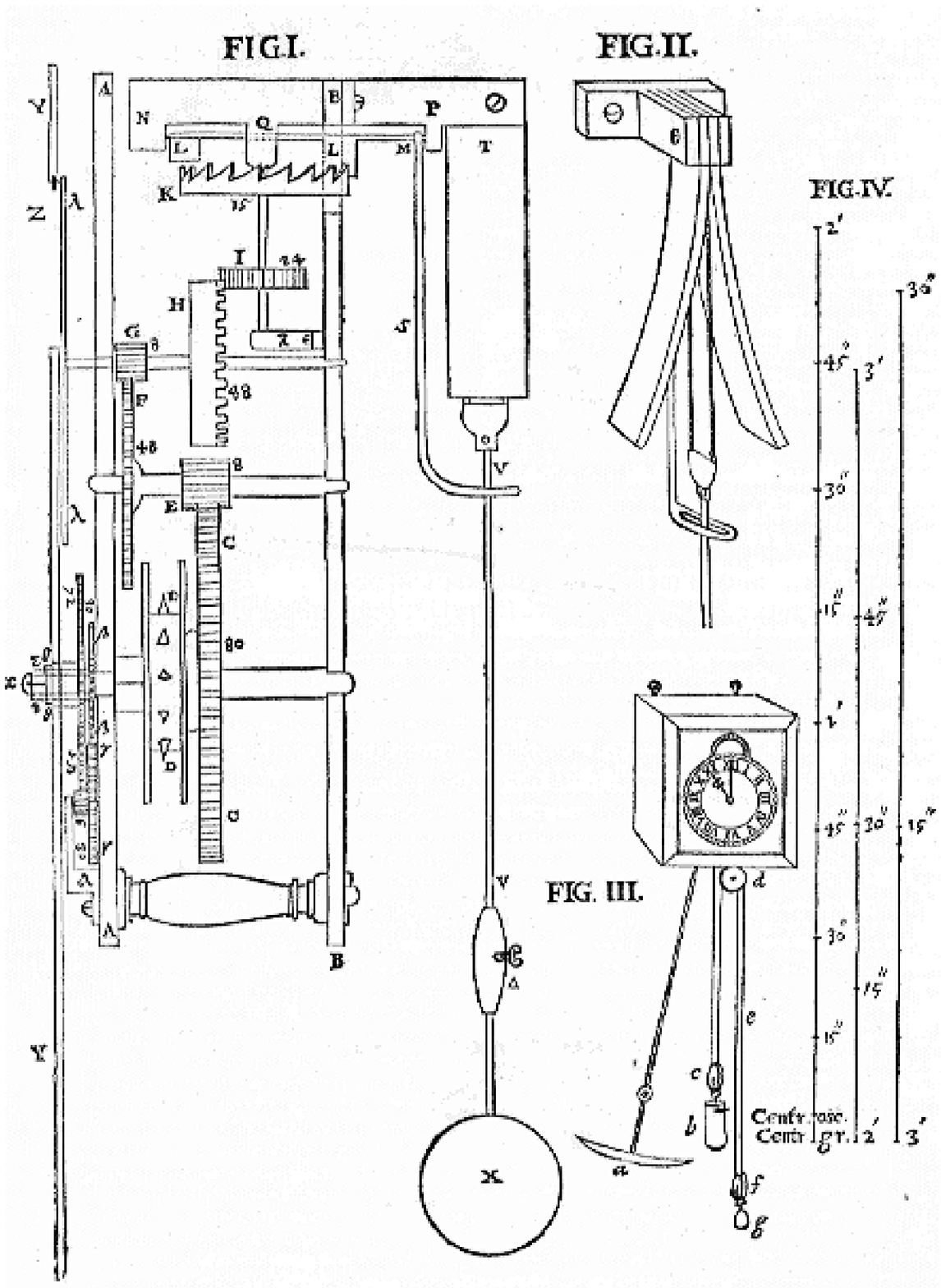}

Figure 5
\end{center}
\end{minipage}
\vspace{2mm}

\onderwerp{Geometric optics}

Either of the following fundamental principles can be used to explain
the propagation of light:

\begin{Fermat}
{\rm (1658)} Any ray of light follows the path of shortest time.
\end{Fermat}

\begin{Huygens}
{\rm (1690, \cite{huyg90})}
Every point of a wave front is the source of an elementary wave. The wave
front at a later time is given as the envelope of these elementary waves.
\end{Huygens}

The simplest possible example is the propagation of light in a homogeneous
and isotropic medium. Here we expect the rays of light to be straight
lines. Figure~6 illustrates that this is indeed what the two
principles predict. We merely need to observe that,
in a homogeneous and isotropic medium, the curves of shortest time are
the same as geometrically shortest curves, i.e.\ straight lines, and
elementary waves are circular waves around their centre.

\vspace{2mm}
\begin{minipage}{8cm}
\begin{center}
\epsfxsize 5cm \epsfbox{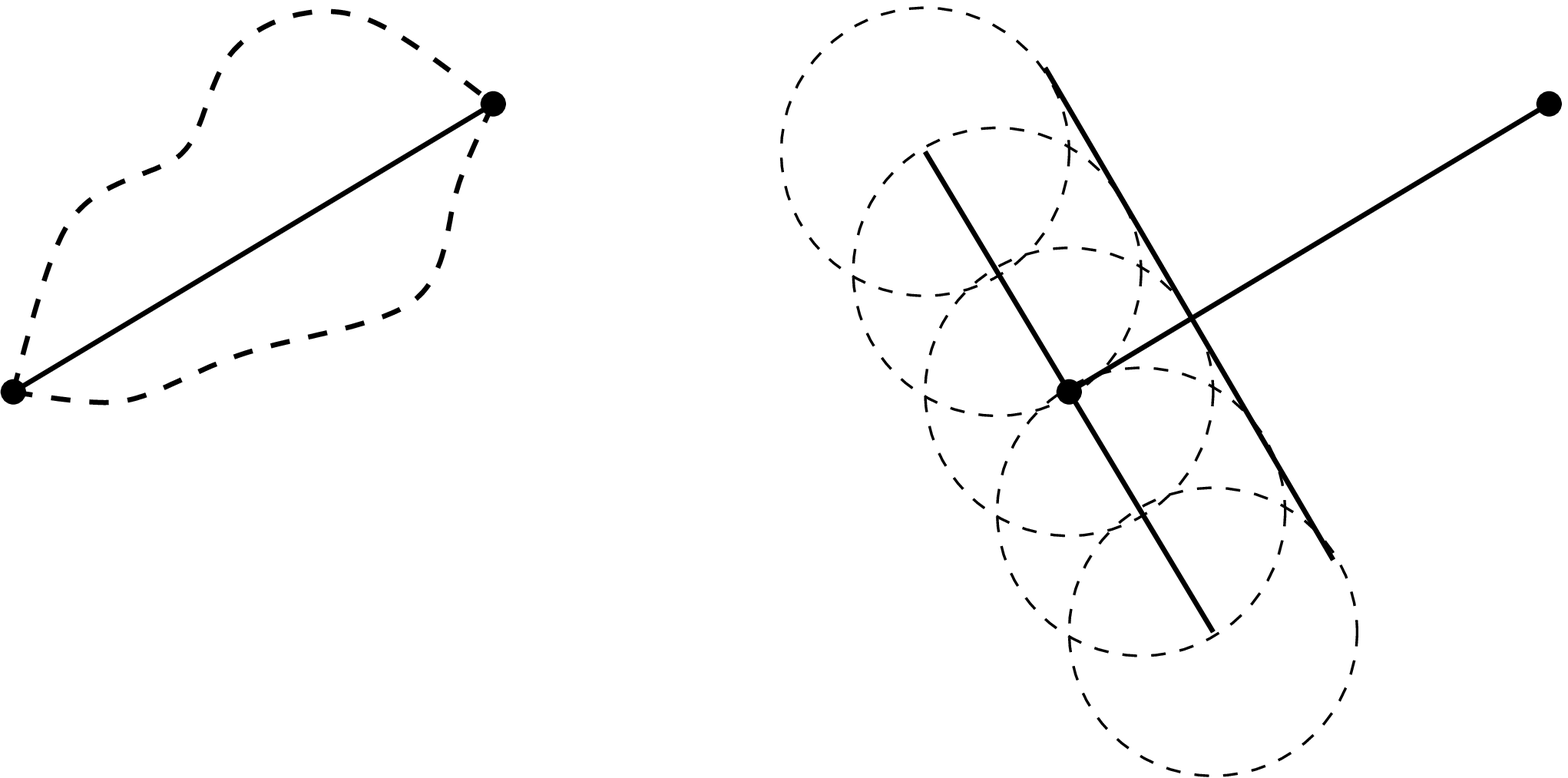}

Figure 6
\end{center}
\end{minipage}
\vspace{2mm}

Whereas Fermat's principle can only be justified as an instance of nature's
parsimony, cf.~\cite{hitr96}, Huygens' principle can be explained
mechanistically from a particle theory of light, see Figure~7.

\vspace{2mm}
\begin{minipage}{8cm}
\begin{center}
\epsfxsize 4cm \epsfbox{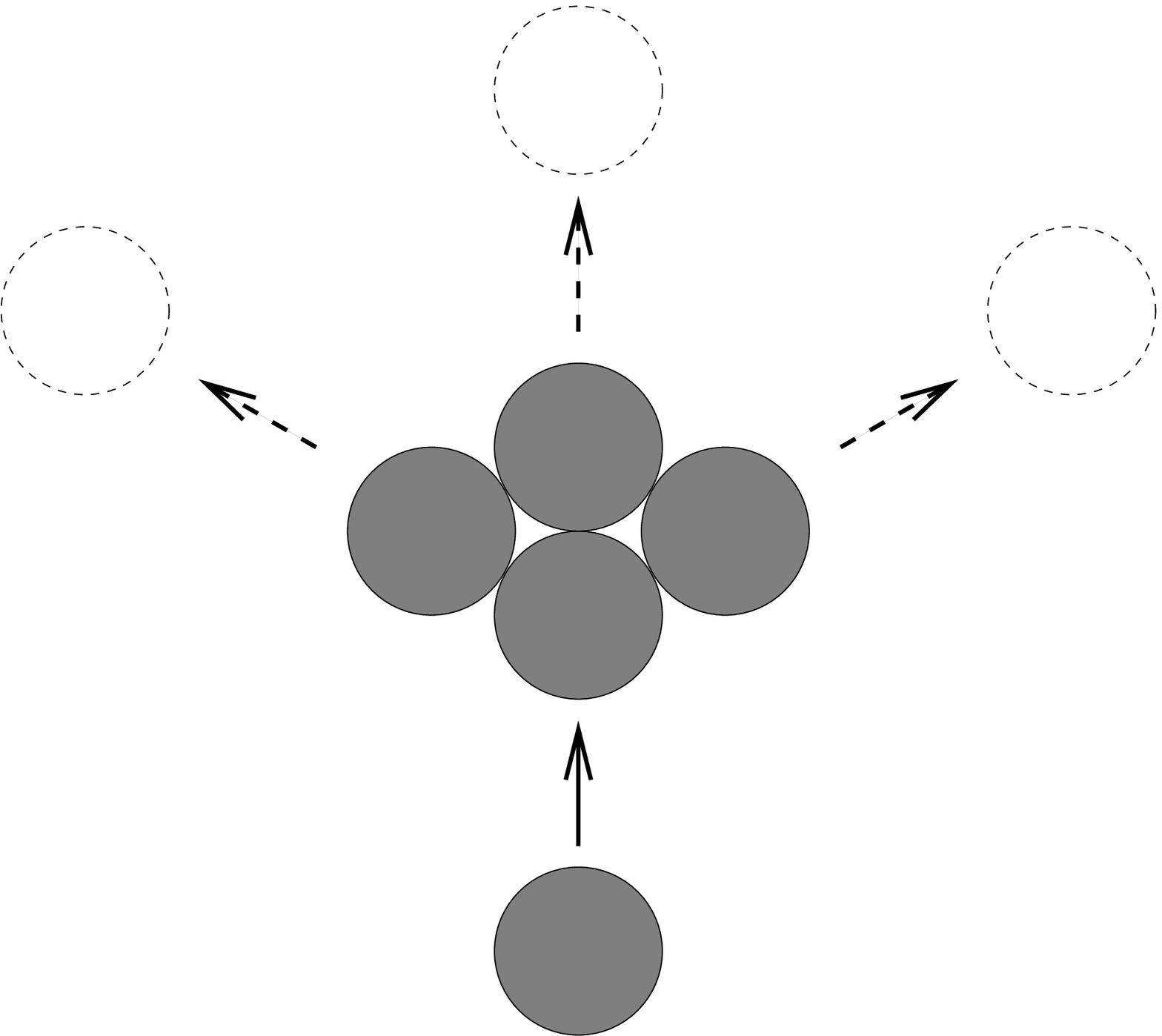}

Figure 7
\end{center}
\end{minipage}
\vspace{2mm}

To illustrate the power of these principles, here are two further examples.
The first is the law of reflection, which states that the angle of
incidence equals the angle of reflection. Figure~8 shows how this follows
from Fermat's principle: The path connecting $A$ and $B$ has the
same length as the corresponding one connecting $A$ and the mirror
image $B'$ of~$B$, and for the latter the shortest (and hence
quickest) path is given by the straight line.

\vspace{2mm}
\begin{minipage}{8cm}
\begin{center}
\centerline{\relabelbox\small
\epsfxsize 8cm \epsfbox{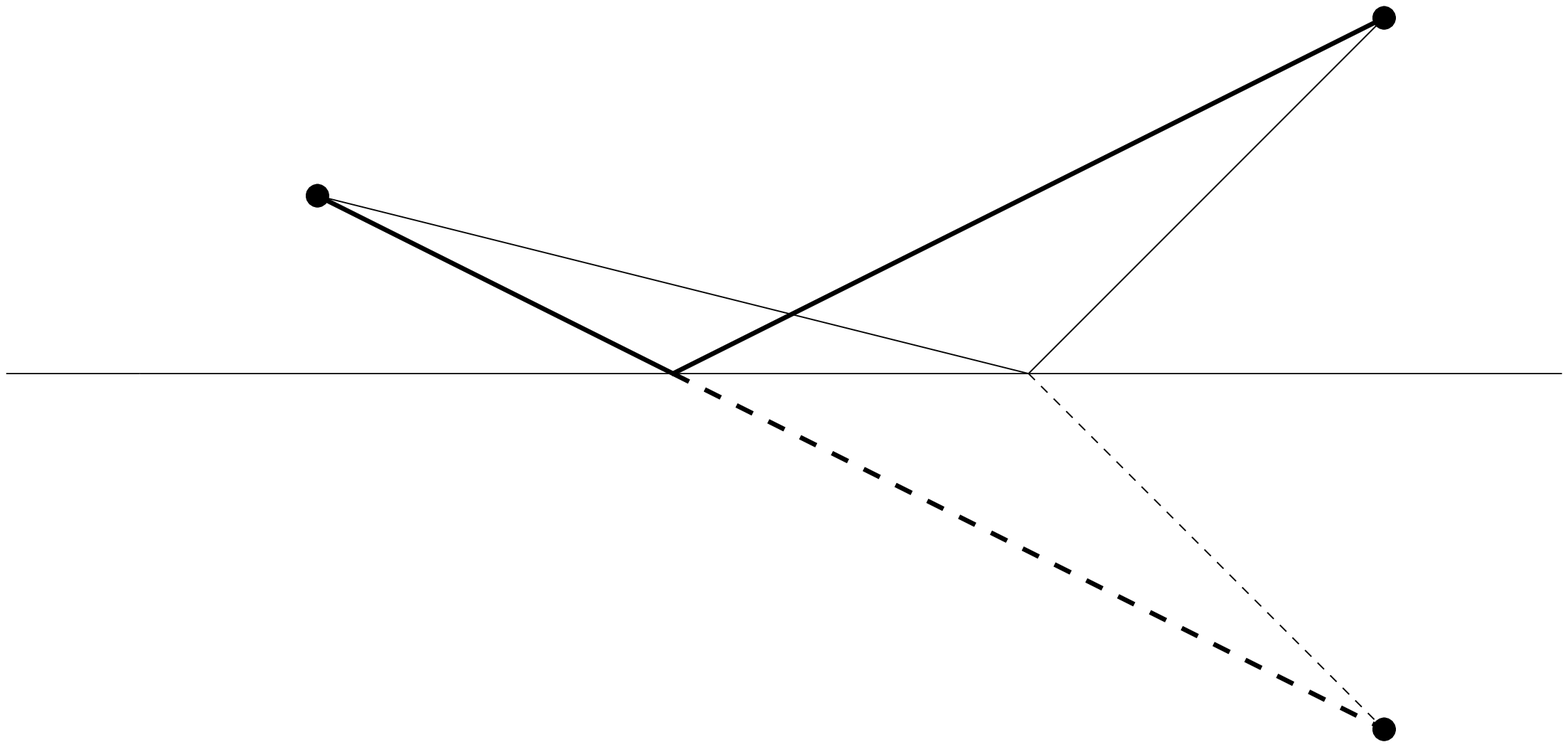}
\extralabel <-1.5cm, 2.0cm> {mirror}
\extralabel <-6.9cm, 2.8cm> {$A$}
\extralabel <-0.8cm, 0.0cm> {$B'$}
\extralabel <-0.8cm, 3.6cm> {$B$}
\endrelabelbox}

Figure 8
\end{center}
\end{minipage}
\vspace{2mm}

The explanation of the law of reflection from Huygens' principle
is illustrated in Figure~9.

\vspace{2mm}
\begin{minipage}{8cm}
\begin{center}
\centerline{\relabelbox\small
\epsfxsize 8cm \epsfbox{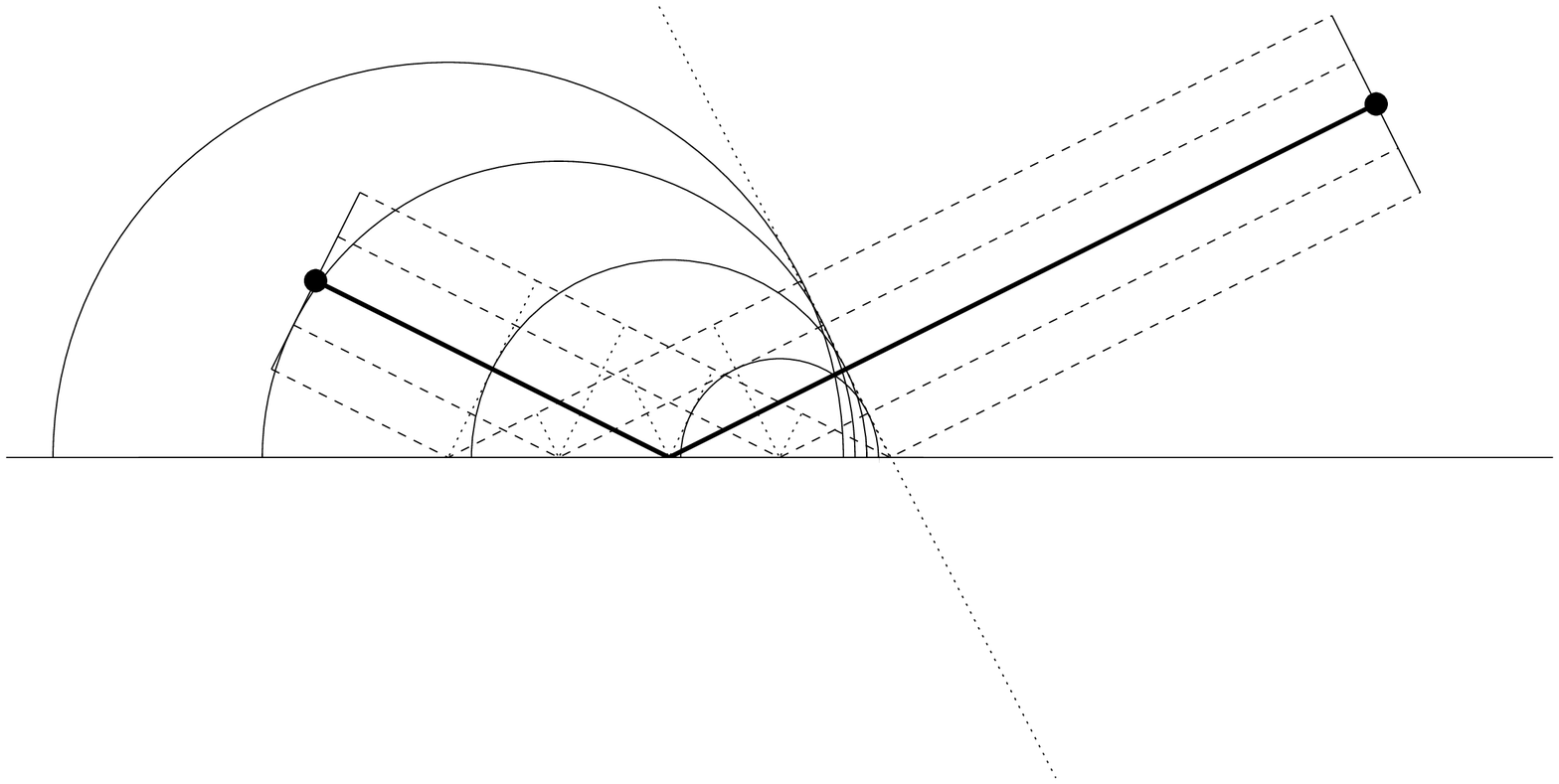}
\extralabel <-1.5cm, 2.4cm> {mirror}
\extralabel <-6.9cm, 3.2cm> {$A$}
\extralabel <-0.8cm, 4.1cm> {$B$}
\endrelabelbox}

Figure 9
\end{center}
\end{minipage}
\vspace{2mm}

As a final application of the two principles, we turn to the law of
refraction, also known as Snell's law after the Dutch astronomer and
mathematician Willebrord van Roijen Snell (1580--1628), whose latinised
name Snellius now adorns the Mathematical Institute of the Universiteit
Leiden. Snell discovered this law in 1621; in print it appears for instance
in Huygens' {\it Trait\'e de la lumi\`ere}, with proofs based on either
of the two principles. The law states that as
a ray of light crosses the boundary between two (homogeneous and
isotropic) optical media, the angle of incidence $\alpha_1$ (measured
relative to a line perpendicular to the separating surface) and the angle
of refraction $\alpha_2$ (see Figure~10) are related by
\[ \frac{\sin\alpha_1}{v_1}=\frac{\sin\alpha_2}{v_2},\]
where $v_1$ and $v_2$ denote the speed of light in the respective medium.

\vspace{2mm}
\begin{minipage}{8cm}
\begin{center}
\centerline{\relabelbox\small
\epsfxsize 8cm \epsfbox{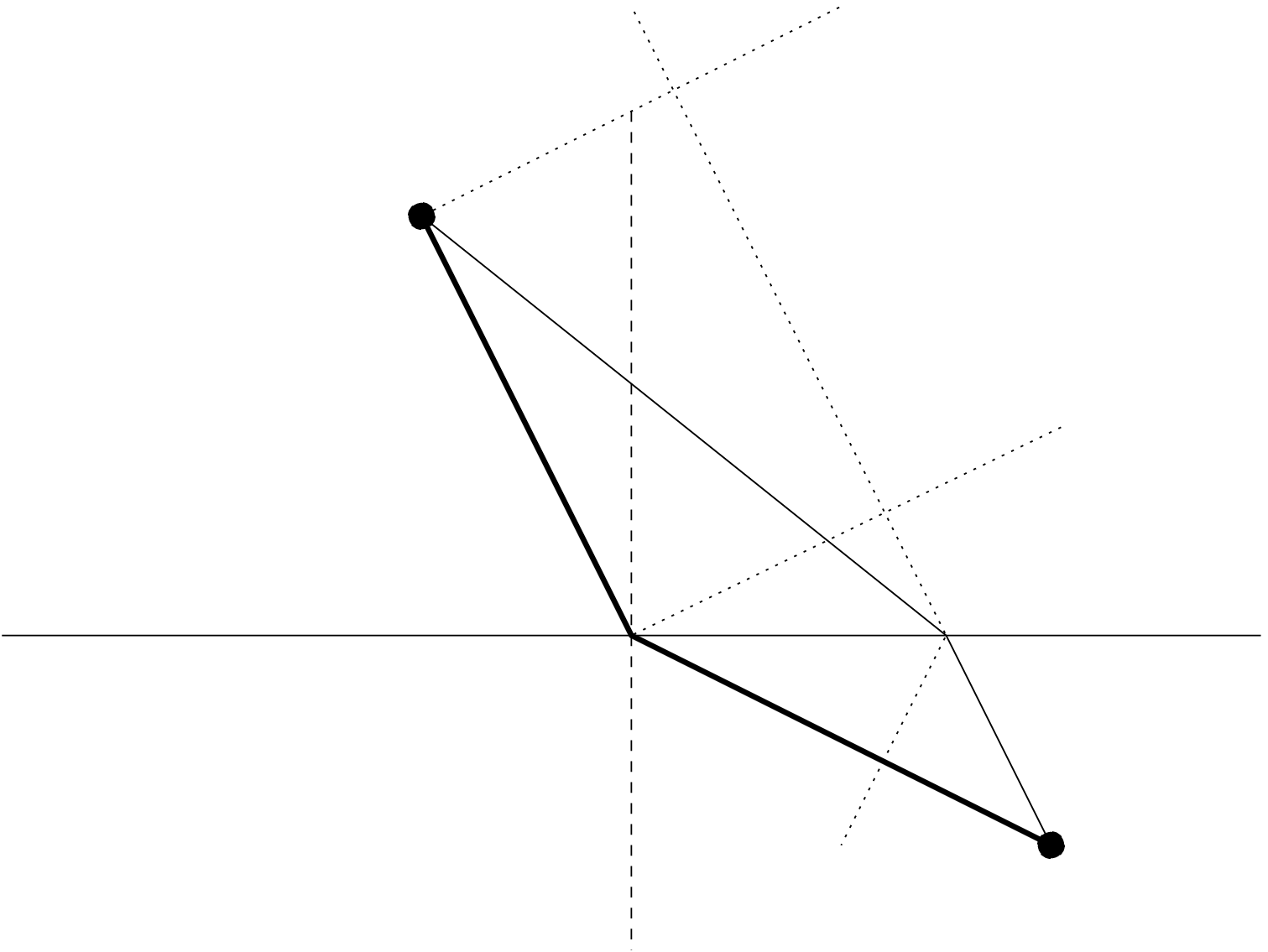}
\extralabel <-5.9cm, 4.6cm> {$A$}
\extralabel <-3.9cm, 5.7cm> {$A'$}
\extralabel <-2.5cm, 3.1cm> {$Q'$}
\extralabel <-2.0cm, 2.2cm> {$Q$}
\extralabel <-4.5cm, 1.7cm> {$P$}
\extralabel <-2.6cm, 0.8cm> {$P'$}
\extralabel <-4.5cm, 3.0cm> {$\scriptstyle \alpha_1$}
\extralabel <-4.0cm, 1.6cm> {$\scriptstyle \alpha_2$}
\extralabel <-1.2cm, 0.6cm> {$B$}
\extralabel <-7.5cm, 2.15cm> {medium 1}
\extralabel <-7.5cm, 1.7cm> {medium 2}
\endrelabelbox}

Figure 10
\end{center}
\end{minipage}
\vspace{2mm}

Figure 10 shows how to derive Snell's law from Fermat's principle. The path
from $A$ to $B$ via $P$ (drawn in bold) is supposed to be the one satisfying
Snell's law. We need to show that it takes longer to travel along any
other broken path from $A$ to $B$ via some $Q$ different from~$P$.
We compute
\[ \frac{PP'}{v_2}=\frac{PQ\sin\alpha_2}{v_2}=\frac{PQ\sin\alpha_1}{v_1}
=\frac{QQ'}{v_1},\]
that is,
\[ t(PP') =t(QQ'),\]
where $t(\cdot )$ denotes the amount of time it takes to travel along a
certain line segment in the corresponding medium. Therefore
\begin{eqnarray*}
t(AQ)+t(QB) & > & t(A'Q)+t(P'B)\\
            & = & t(A'Q')+t(Q'Q)+t(P'B)\\
            & = & t(AP)+t(PP')+t(P'B)\\
            & = & t(AP)+t(PB).
\end{eqnarray*}

Figure 11 indicates how Snell's law is implied by Huygens' principle.

\vspace{2mm}
\begin{minipage}{8cm}
\begin{center}
\centerline{\relabelbox\small
\epsfxsize 8cm \epsfbox{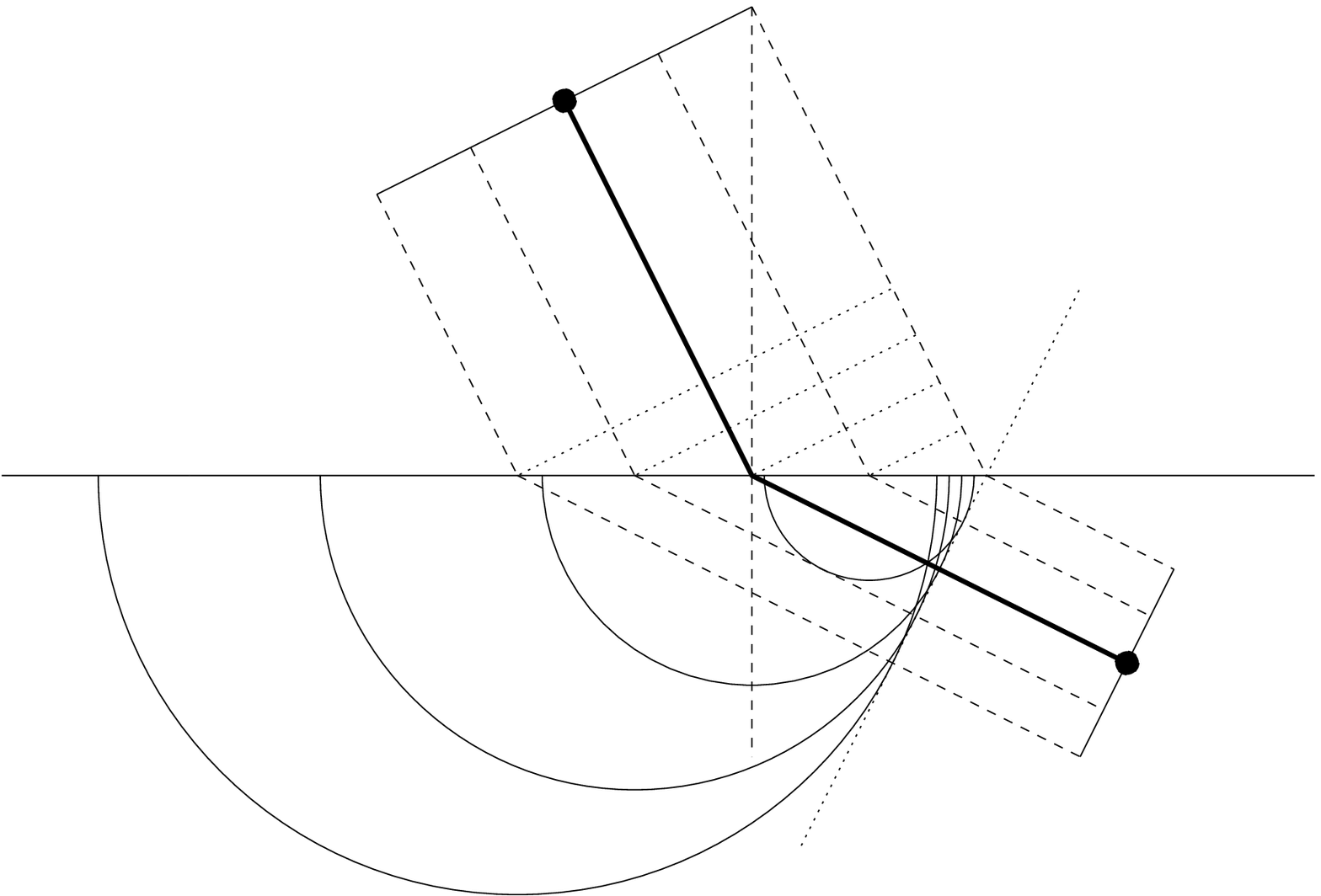}
\extralabel <-1.5cm, 2.7cm> {medium 1}
\extralabel <-1.5cm, 2.25cm> {medium 2}
\endrelabelbox}

Figure 11
\end{center}
\end{minipage}
\vspace{2mm}

\onderwerp{Johann Bernoulli's solution of the brachistochrone problem}

Jacob Bernoulli, in his response (cited in Section~2) to the fraternal
challenge, developed a general method for dealing with problems of this kind,
nowadays known as the calculus of variations. In the present section we shall
be concerned with Johann's own solution, which nicely relates to the concepts
of geometric optics discussed above.

When the mass $M$ has reached a point $(x,y)$ on the slide from $A=(0,0)$
to $B$, with the $y$--coordinate oriented downwards, its speed has reached,
under the influence of gravitation, the value
\[ v=\sqrt{2gy},\]
where $g=9.81m/s^2$ denotes the gravitational acceleration near the
surface of the Earth. In order to determine which path the point $M$ should
follow so as to take the shortest time from $A$ to~$B$, we discretise
the problem.

Imagine that the region between $A$ and $B$ is layered into finitely
many horizontal slices, in each of which the speed of $M$ stays constant.
In particular, $M$ should follow a straight line in each layer. As $M$
passes from the $i$th to the $(i+1)$st layer, the angle $\alpha_i$ of
incidence and $\alpha_{i+1}$ of `refraction' should be related to the
respective speeds $v_i,v_{i+1}$ by Snell's law
\[ \frac{v_{i+1}}{\sin\alpha_{i+1}}=\frac{v_i}{\sin\alpha_i},\]
for the fact that Snell's law is an instance of Fermat's principle
guarantees this to yield the quickest path (Figure~12).

\vspace{2mm}
\begin{minipage}{8cm}
\begin{center}
\centerline{\relabelbox\small
\epsfxsize 8cm \epsfbox{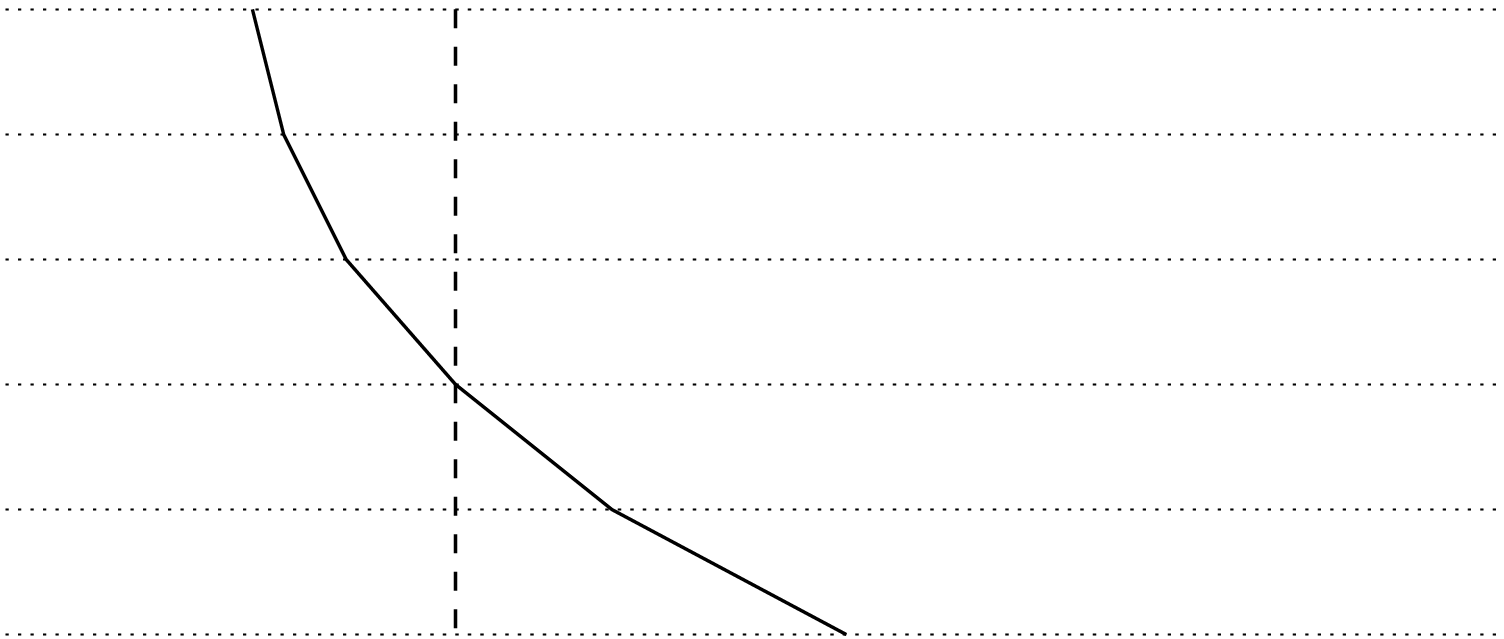}
\extralabel <-1.5cm, 1.6cm> {$v_i$}
\extralabel <-1.5cm, 0.9cm> {$v_{i+1}$}
\extralabel <-6.0cm, 1.8cm> {$\scriptstyle \alpha_i$}
\extralabel <-5.6cm, 0.8cm> {$\scriptstyle \alpha_{i+1}$}
\endrelabelbox}

Figure 12
\end{center}
\end{minipage}
\vspace{2mm}

As we let the number of slices
tend to infinity, the equation describing the brachistochrone becomes
\[ \frac{v}{\sin\alpha}=c\]
for some constant $c$, see Figure~13. Bravely computing with infinitesimals,
we have $\sin\alpha = dx/\sqrt{dx^2+dy^2}$, whence
\[ \sqrt{1+\bigl(\frac{dy}{dx}\bigr)^2}\cdot\sqrt{2gy}=c.\]
This can be written as
\[ \frac{dx}{dy} = \sqrt{\frac{y}{2a-y}}\]
with $a=c^2/4g$. Substitute
\[ y(t)=2a\sin^2\frac{t}{2}=a (1-\cos t).\]
Then
\[ \frac{dx}{dt}=\frac{dx}{dy}\cdot\frac{dy}{dt}=
\sqrt{\frac{1-\cos t}{1+\cos t}}\cdot a\sin t=2a\sin^2\frac{t}{2}=y,\]
hence (with $x(0)=0$)
\[ x(t)=a(t-\sin t),\]
``ex qua concludo: curvam {\it Brachystochronam} esse {\it Cycloidem
vulgarem}''~\cite[p.~266]{bern}.

\vspace{2mm}
\begin{minipage}{8cm}
\begin{center}
\centerline{\relabelbox\small
\epsfxsize 6cm \epsfbox{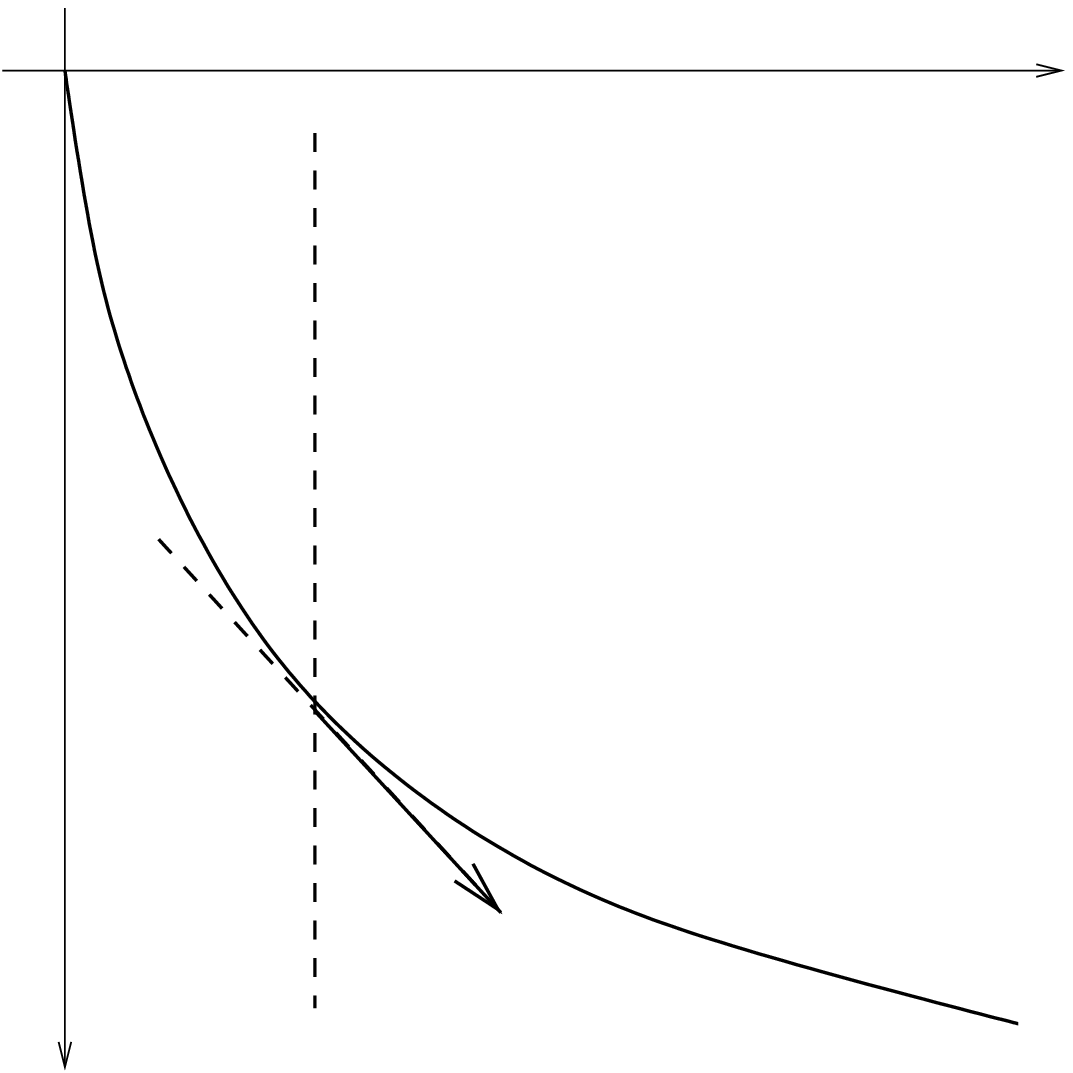}
\extralabel <-0.3cm, 5.8cm> {$x$}
\extralabel <-6.1cm, 0.1cm> {$y$}
\extralabel <-4.2cm, 1.3cm> {$\alpha$}
\extralabel <-3.4cm, 0.7cm> {$v$}
\endrelabelbox}

Figure 13
\end{center}
\end{minipage}
\vspace{2mm}

This is as good a point as any to recommend the wonderful
textbook~\cite{hawa00}. It contains an extensive discussion of both
the brachistochrone and tautochrone problem in their historical
context, and many other historical gems that so sadly are missing
from our usual introductory courses on analysis, which tend to suffer from
the dictate of efficiency and the haste to `cover material'.

\onderwerp{Elementary contact geometry}

Here at last we come to the second part of this article's title.
My modest aim is to convey a couple of basic notions of contact geometry
and to show how they relate to some of the ideas discussed above. In doing
so, I am aware of W.~Thurston's warning that ``one person's clear
mental image is another person's intimidation''~\cite{thur94}.

One of the fundamental notions of contact geometry is the so-called
{\it space of (oriented) contact elements} of a given manifold. Let us
first consider a concrete example, see Figure~14.

\vspace{2mm}
\begin{minipage}{8cm}
\begin{center}
\centerline{\relabelbox\small
\epsfxsize 8cm \epsfbox{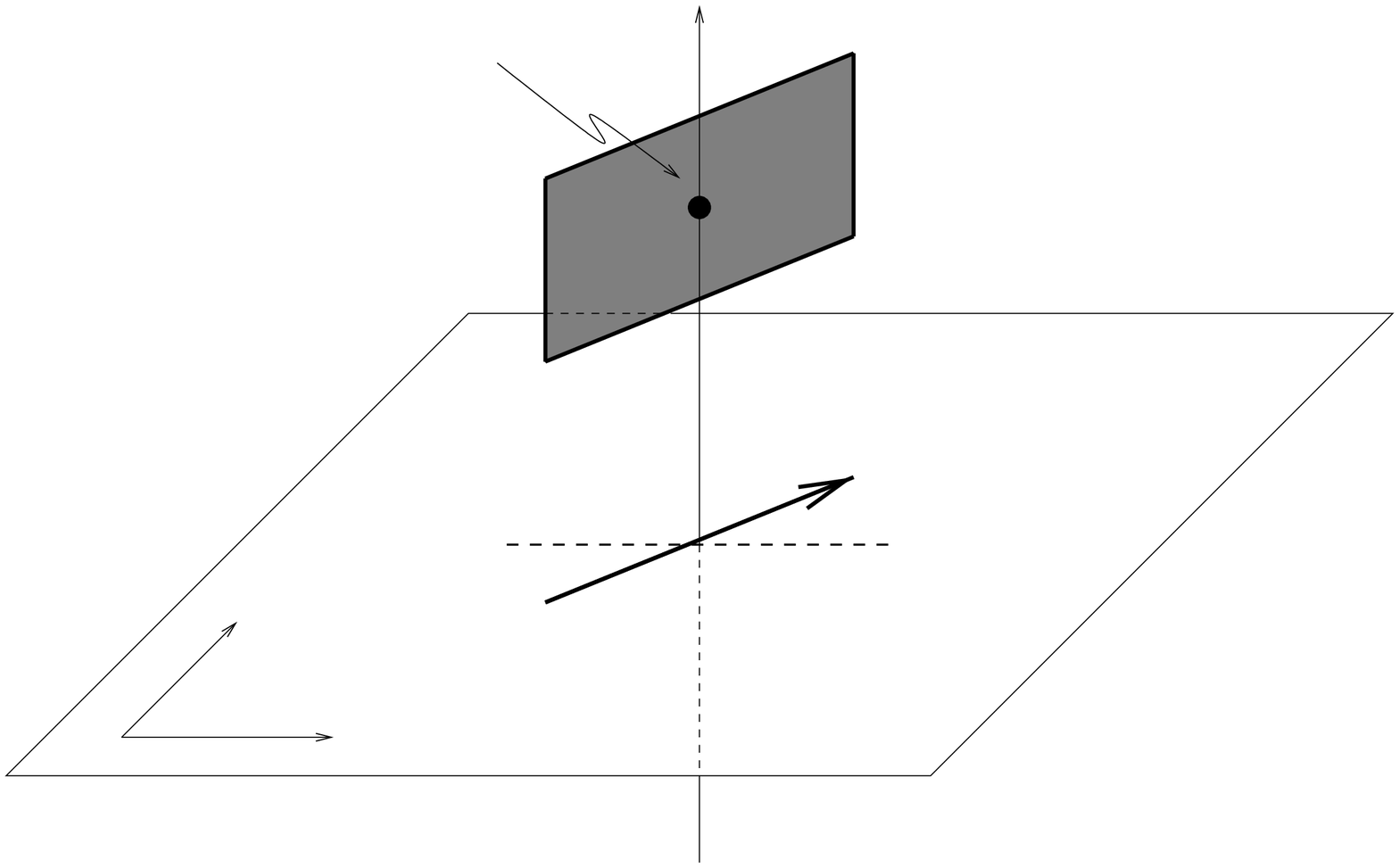}
\extralabel <-6.1cm, 0.7cm> {$x$}
\extralabel <-6.65cm, 1.3cm> {$y$}
\extralabel <-4.0cm, 1.5cm> {$p_0$}
\extralabel <-3.4cm, 1.95cm> {$\scriptstyle\theta_0$}
\extralabel <-4.4cm, 4.8cm> {$\theta$}
\extralabel <-5.8cm, 4.7cm> {$(p_0,\theta_0)$}
\extralabel <-3.0cm, 4.0cm> {$\xi_{(p_0,\theta_0)}$}
\endrelabelbox}

Figure 14
\end{center}
\end{minipage}
\vspace{2mm}

An oriented contact element to the $2$--plane
$\R^2$ at some point $p_0\in\R^2$ is simply an oriented line
passing through the point~$p_0$. Such a line is uniquely determined
by $p_0$ and an angle~$\theta_0$.
We can think of this
angle $\theta_0$ as an element of the unit circle~$S^1$, so the space of
all contact elements of $\R^2$ can be identified with the product
$\R^2\times S^1$.

Let $\partial_x,\partial_y$ denote
the unit vectors in the coordinate directions of~$\R^2$ at any
given point $(x_0,y_0,\theta_0)\in\R^2\times S^1$. They can be thought
of as the velocity vectors of the curves
\[ x\longmapsto (x_0+x,y_0,\theta_0) \;\;\; \mbox{\rm resp.}\;\;\;
y\longmapsto (x_0,y_0+y,\theta_0).\]
Similarly, we can speak of the velocity vector $\partial_{\theta}$
of the curve
\[ \theta\longmapsto (x_0,y_0,\theta_0+\theta ).\]
In the local picture of Figure~14, where $\theta$ is measured along a real
axis, this is once again simply the unit vector in the direction of
the $\theta$--coordinate.

We now specify a $2$--plane $\xi_{(p_0,\theta_0)}$ at any point $(p_0,\theta_0)
\in\R^2\times S^1$ as the plane spanned by the vectors
\[ \partial_{\theta}\;\;\;\mbox{\rm and}\;\;\; \cos\theta_0\,\partial_x
+\sin\theta_0\,\partial_y.\]
Alternatively, this plane is determined by the condition that it
contain $\partial_{\theta}$ and that it project to the contact element
at $p_0$ defined by~$\theta_0$. The collection of all these $2$--planes
is called the {\it natural contact structure} on the space
of contact elements of~$\R^2$.

This probably sounds esoteric or banal, depending on your education.
It is anavoidable that at this point I shall have to assume a certain
level of mathematical literacy. But I make no apology for continuing,
whenever possible, also to address those who are beginning to
feel just a little intimidated. The space of contact elements
of $\R^2$ in fact has a very natural interpretation as a space
of physical configurations. If you want to
describe the position of a wheel of your bicycle, say, you
should describe its position $p$ and its direction, given by~$\theta$.
Moreover, instantaneously the wheel can only roll in the direction in
which it points at any given moment, so the motion of the wheel, interpreted as
a curve in the $3$--dimensional space of contact elements, will be
tangent to the natural contact structure. 

The next concept we want to introduce is that of a
{\it contact transformation}. Such transformations play an important role in
the geometric theory of differential equations. Most physicists first
encounter them in their special incarnation as so-called
Legendre transformations. For our purposes, we can define a contact
transformation as a diffeomorphism $\phi$ of the space of contact elements
$\R^2\times S^1$ with the property that if a curve $w$ passes through
a point $(p,\theta )$ and is tangent to the $2$--plane $\xi_{(p,\theta )}$
at that point, then the image curve $\phi\circ w$ will be tangent to
$\xi_{\phi (p,\theta )}$ at $\phi (p,\theta )$.

Here is an example
of a whole family of contact transformations:
For $t\in\R$, define
\[ \begin{array}{rrcl}
\phi_t: & \R^2\times S^1 & \longrightarrow & \R^2\times S^1\\
        & (x,y,\theta )  & \longmapsto     & (x-t\sin\theta ,y+t\cos\theta ,
                                             \theta ).
\end{array}\]
In order to verify that these are indeed contact transformations, consider
a parametrised curve
\[ s\longmapsto w(s)=(x(s),y(s),\theta (s))\in\R^2\times S^1,\;\;
s\in (-\varepsilon ,\varepsilon ), \]
for some small $\varepsilon >0$ say, with tangent vector
\[ w'(0)=(x'(0),y'(0),\theta '(0))\]
assumed to lie in $\xi_{w(0)}$.
With $\pi\colon\,\R^2\times S^1\rightarrow\R^2$ denoting the natural
projection, this is equivalent to
saying that the tangent vector $(x'(0),y'(0))$ of the projected curve
$\pi\circ w$ at the point $(x(0),y(0))$
lies in the line determined by~$\theta (0)$, i.e.\ is a multiple of
$(\cos\theta (0),\sin\theta (0))$.

The transformed curve is
\[ \phi_t\circ w(s)=(x(s)-t\sin\theta (s), y(s)+t\cos\theta (s),\theta (s)).\]
Notice that the $\theta$--coordinate remains unchanged under~$\phi_t$.
We compute
\begin{eqnarray*}
\lefteqn{\frac{d}{ds}(\phi_t\circ w)(s)=} \\
 & & (x'(s)-t\theta '(s)\cos\theta (s),
y'(s)-t\theta '(s)\sin\theta (s),\theta '(s))
\end{eqnarray*}
and observe that the $\R^2$--component of this vector at $s=0$
does again lie in the line determined by~$\theta (0)$.

This family $\phi_t$ of transformations is called the {\it geodesic flow}
of~$\R^2$. Here is why: In a general Riemannian manifold, geodesics are
locally shortest curves. In $\R^2$ (with its euclidean metric),
therefore, geodesics are simply the straight lines. Given a point $p\in\R^2$
and a direction $\theta\in S^1$ defining a contact element, let
$\ell_{p,\theta}$ be the unique oriented line in $\R^2$ passing
through the point $p$ and positively orthogonal to the contact
element~$\theta$. This
line is parametrised by
\[ t\longmapsto p+t(-\sin\theta ,\cos\theta ),\;\; t\in\R .\]
Lo and behold, this is the same as $t\mapsto\pi\circ\phi_t(p,\theta )$.
The $\theta$--component of $\phi_t(p,\theta )$ encodes the direction orthogonal
to this geodesic; in our case this component stays constant.

Great, I hear you say, but what does all that have to do with Huygens? Well,
it turns out that we are but one simple step away from proving, with the
help of contact geometry, the equivalence of the principles of Fermat
and Huygens.

Let $\overline{f}$ be a wave front in~$\R^2$, thought of as a parametrised
curve $s\mapsto (x(s),y(s)$, $s\in (-\varepsilon ,\varepsilon )$.
For simplicity, we assume this to be regular, i.e.\
\[\overline{f}'(s)=(x'(s),y'(s))\neq (0,0)\;\;\mbox{\rm for all} \;\;
s\in (-\varepsilon ,\varepsilon ).\]
Such a wave front lifts to a unique curve
\[ s\longmapsto f(s)=(x(s),y(s),\theta (s)) \]
in the space of contact elements subject to the requirement that
$(x'(s),y'(s))$ be a positive multiple of $(\cos\theta (s),\sin\theta (s))$;
this lift will be tangent to the natural contact structure.
{\bf Fermat's principle} says that light propagates along the geodesic
rays (i.e.\ straight lines) orthogonal to the wave front~$\overline{f}$,
which translates into saying that the wave front at some later time $t$ is
given by $\pi\circ\phi_t\circ f$.

Next consider the curve
\[ h\colon\,\theta\longmapsto (x(0),y(0),\theta ).\]
This is
simply the circle worth of all contact elements at the point
$\pi\circ h\equiv (x(0),y(0))$. Under the geodesic flow and projected
to~$\R^2$, this becomes an elementary wave in the sense of Huygens:
for each fixed $t\in\R$ the curve
\[ \theta\longmapsto\pi\circ\phi_t\circ h(\theta )=
(x(0),y(0))+t(-\sin\theta ,\cos\theta )\]
is a circle of radius~$t$ centred at $(x(0),y(0)$.

The curves $h$ and $f$ are both tangent to $\xi_{f(0)}$ at the point
$f(0)=h(\theta (0))$. Since $\phi_t$ is a contact transformation,
the transformed curves $\phi_t\circ h$ and $\phi_t\circ f$ will
be tangent to $\xi_{\phi_t\circ f(0)}$ at $\phi_t\circ f(0)$. Then, by the
definition of the natural contact structure, the transformed wave front
$\pi\circ\phi_t\circ f$ and the elementary wave $\pi\circ\phi_t\circ h$
will be tangent to each other at the point $\pi\circ\phi_t\circ f(0)$ ---
this is {\bf Huygens' principle}.

The general argument is entirely analogous: A contact element on a
Riemannian manifold is a (cooriented) tangent hyperplane field. The space of
all these contact elements once again carries a natural contact
structure. A geodesic is uniquely determined by an initial point
and a direction positively orthogonal to a contact element at that
point. Like in the special case of $\R^2$
one can show that the geodesic flow preserves the natural contact structure
on the space of contact elements, and this translates into the equivalence
of the two principles of geometric optics. A quick proof of this
general case is given in~\cite{geig01}; full details of that proof
are meant to appear in a forthcoming book on contact topology.

\onderwerp{The family portrait}

It remains to identify the young Christiaan Huygens in Hanneman's family
portrait. In the biography~\cite{andr94} (from an aptly named
publishing company!), a whole chapter is devoted
to this question, so we seem to be in muddy waters.

Since Christiaan was the second-eldest son, there is actually only a choice
between the two boys at the top. My first guess was that Christiaan is
the one on the left, who has arguably the most striking face. This intuitive
feeling is confirmed by the catalogue of the
Mauritshuis~\cite[p.~67]{maur93} and by the afterword in~\cite{huyg90}.
Alas, it is wrong.

It appears that the confusion was started by an engraving of the printing
carried out for a late 19th century edition of the collected works
of Christiaan Huygens. Here Christiaan's name is placed at the upper left,
contradicting an earlier engraving; the original painting does not
associate names with the four boys. However, family iconography of the time
demanded that the eldest son be placed to his father's right, i.e.\ on the
left side of the portrait. This identification of the eldest brother
Constantijn as the boy on the upper left, and thus Christiaan as the
one on the right, seems to be confirmed by a comparison of the painting
with other portraits from the same period.

\completepublications

\end{document}